\documentclass{elsart}
\def\Z/{{\mathbb Z}}
\def\Zn/{{\mathbb Z}_n}
\def\beh/{behaving}
\def\zfs/{zero-free sequence}
\def\zf/{zero-free}
\def\pin/{positive integer}
\def\pins/{positive integer sequence}
\def\lpr/{least positive representatives}
\def\v#1{\overline{#1}}

\def\sumb{\sum_{j=1}^q\v{-b_j}}
\def\intpart{\left\lfloor\left(\frac{k-1}{2}\right)^2\right\rfloor}
\usepackage{amssymb,amsmath}

\begin{document}
\begin{frontmatter}
\title{Long ${\mathbf n}$-\zfs/s in finite cyclic groups}
\author[]{Svetoslav Savchev}
\kern-.35em\footnote{No current affiliation.},
\author[Emory]{Fang Chen}
\address[Emory]{Oxford College of Emory University, Oxford, GA 30054, USA}

\begin{abstract}
A sequence in the additive group~$\Zn/$ of integers modulo~$n$ is
called $n$-\zf/ if it does not contain subsequences with
length~$n$ and sum zero. The article characterizes the $n$-\zfs/s
in~$\Zn/$ of length greater than~$3n/2{-}1$. The structure of
these sequences is completely determined, which generalizes a
number of previously known facts. The characterization cannot be
extended in the same form to shorter sequence lengths.
Consequences of the main result are best possible lower bounds for
the maximum multiplicity of a term in an $n$-\zfs/ of any given
length greater than~$3n/2{-}1$ in~$\Zn/$, and also for the
combined multiplicity of the two most repeated terms. Yet another
application is finding the values in a certain range of a function
related to the classic theorem of Erd\H{o}s, Ginzburg and~Ziv.
\end{abstract}

\end{frontmatter}

\section{Introduction}
\label{Intro}

The Erd\H{o}s--Ginzburg--Ziv theorem \cite{ErdosGinzburgZiv}
states that each sequence of length ${2n{-}1}$ in the cyclic group
of order~$n$ has a subsequence of length~$n$ and sum~zero. This
article characterizes all sequences with length greater
than~$3n/2{-}1$ in the same group that do not satisfy the
conclusion of the celebrated theorem.

In the sequel, the cyclic group of order~$n$ is identified with the
additive group~$\Zn/=\Z//n\Z/$ of integers modulo~$n$.
A sequence in~$\Zn/$ is called a {\em
zero sequence} or a {\em zero sum} if its sum is the zero element
of~$\Zn/$. A sequence is {\em \zf/} if it does not contain a
nonempty zero subsequence. Sequences in~$\Zn/$ without zero
subsequences of length~$n$ will be called {\em $n$-\zf/}.

\looseness=-1 The $n$-\zfs/s in~$\Zn/$ were given considerable
attention. Here we mention only results most closely related to
our topic. Yuster and~Peterson \cite{PetersonYuster} and,
independently, Bialostocki and Dierker~\cite{BialostockiDierker},
determined all $n$-\zfs/s of length~$2n{-}2$ in~$\Zn/$. These are
the sequences containing exactly two distinct elements $a$ and~$b$
of~$\Zn/$, each of them repeated $n{-}1$~times, such that ${a-b}$
generates~$\Zn/$. Ordaz and Flores~\cite{FloresOrdaz} solved the
same problem for length~$2n{-}3$. Again, two distinct terms have
high combined multiplicity (details can be found in
Section~\ref{Highmult}). In general, the combined multiplicity of
the two most represented terms was intensively studied. Gao
\cite{Gao2} proved a statement to this effect for $n$-\zfs/s of
length roughly greater than~$7n/4$. A recent work of Gao,
Panigrahi and Thangadurai \cite{GaoPanigrahiThagadurai} considered
the same question in the case of a prime~$n$, for length roughly
greater than~$5n/3$.

Based on the main theorem in~\cite{Gao2}, Bialostocki et.~al.\
\cite{BialostockiDierkerGrynkiewiczLotspeich} obtained an explicit
characterization of the $n$-\zfs/s in~$\Zn/$ with length greater
than or equal to ${2n-2-\lfloor n/4\rfloor}$. The core of their
proof is essentially present already in the
article~\cite{GaoHamidoune} of Gao and Hamidoune.

Our goal is to characterize the $n$-\zfs/s in~$\Zn/$ of length
greater than~$3n/2{-}1$. The argument relies on a key structural
result from~\cite{PCheetahSCat} about \zfs/s of length greater
than~$n/2$ in~$\Zn/$. The description obtained generalizes the one
from~\cite{BialostockiDierkerGrynkiewiczLotspeich} and cannot be
extended in the same shape to shorter sequences. In this sense the
range of the characterization is optimal.

Let $a$ be an integer coprime to~$n$ and $b$ an element of~$\Zn/$.
The function $x\mapsto ax+b$ from~$\Zn/$ into itself will be called
an {\em affine map}. In particular the {\em translations}~$x\mapsto
x+b$ are affine maps, for each~$b\in \Zn/$. Suppose that a
sequence~$\beta$ in~$\Zn/$ can be obtained from a sequence~$\alpha$
through an affine map and rearrangement of terms. Then we say that
$\alpha$ is {\em similar} to~$\beta$ and write $\alpha\sim\beta$.
Clearly $\sim$ is an equivalence relation. Affine maps preserve zero
sums of length~$n$ and do not bring in new ones. So it is usual not
to distinguish between similar sequences in questions involving
$n$-term zero subsequences. Our characterization will be up to
similitude, i.~e. up to affine maps and rearrangement of terms.

If $a\in\Zn/$, let $\v{a}$ denote the unique integer in~$[1,n]$ that
belongs to the congruence class~$a$ modulo~$n$. We call~$\v{a}$ the
{\em least positive representative} of~$a$. Least positive
representatives occur frequently in the text, so we allow a certain
abuse of notation to simplify the exposition. In some cases we write
$a$ instead of~$\v{a}$. This applies mostly to the group elements
$0$ and $1$; the actual meaning of the symbols~$0$ and $1$ should be
clear from the context. Furthermore, for a sequence~$\alpha$
in~$\Zn/$ we denote by $\v{\alpha}$ the sequence of its \lpr/, and
by $L(\alpha)$ the sum of~$\v{\alpha}$.

The sequences considered are written multiplicatively, and
multiplicities of sequence terms are indicated by using exponents.
The length of the sequence~$\alpha$ is denoted by~$|\alpha|$. The
{\em union} of two sequences $\alpha$ and $\beta$, denoted
$\alpha\cup\beta$, is the sequence formed by appending the terms
of~$\beta$ to~$\alpha$. Also, $1-\beta$ is the sequence obtained
by replacing each term $b$ of $\beta$ by $1-b$.

Now the main result in the article, Theorem~\ref{characterization},
can be stated as follows:
\begin{quote}\em
A sequence of length greater than $3n/2-1$ in~$\Zn/$ is $n$-\zf/ if
and only if it is similar to the union of two sequences $\alpha$ and
$\beta$ in~$\Zn/$ such that $L(\alpha)<n$ and $L(1-\beta)<n$.
\end{quote}
Once such a characterization is available, certain basic questions
about sufficiently long $n$-\zfs/s in cyclic groups receive
satisfactory answers.

The preliminaries needed for the key proof are included in
Section~\ref{Pre}. The main result is proven in
Section~\ref{MainResult}. It is preceded by some properties of
sequences of the form $\alpha\cup\beta$, where $\alpha$ and $\beta$
are sequences in~$\Zn/$ satisfying $L(\alpha)<n$ and $L(1-\beta)<n$.
In Section~\ref{Highmult} we study the maximum multiplicity of a
term in an $n$-\zfs/ of length~$n-1+k$, where $n/2<k<n$, and also
the combined multiplicity of the two most repeated terms. Best
possible lower bounds are established in both cases. The main
theorem enables us to determine, in Section~\ref{gnk}, the values in
a certain range of a function related to a variant of the
Erd\H{o}s--Ginzburg--Ziv theorem.

\section{Preliminaries}
\label{Pre}

For sequences $\alpha$ and $\beta$ in~$\Zn/$, we say that $\alpha$
is {\em equivalent} to~$\beta$ if $\beta$ can be obtained
from~$\alpha$ through multiplication by an integer coprime to~$n$
and rearrangement of terms. Such multiplication is an affine map
preserving all zero sums in~$\Zn/$, not just the ones of
length~$n$. In particular equivalent sequences are similar. Our
characterization rests on the following result
from~\cite{PCheetahSCat}:
\begin{thm}\label{zero-free}
Each \zfs/ of length greater than~$n/2$ in the cyclic group~$\Zn/$
is equivalent to a sequence whose sum of the \lpr/ is less than~$n$.
\end{thm}

A restatement of a fact from~\cite{Gao1} is also used in the main
proof.

\begin{prop}\label{gao}
A sequence in an abelian group of order~$n$ is such that the
multiplicity of each term is at most the multiplicity
of\kern1pt~$0$. Then each subsequence sum of length greater
than~${n}$ equals a subsequence sum of length exactly~${n}$.
\end{prop}

One more statement is necessary for the main argument.

\begin{prop}\label{general}
Let $\alpha$ be a sequence with positive integer terms of
length~$\ell$ and sum~$S$, where ${2\ell>S}$. Then:

a) $\alpha$ contains at least $2\ell-S$ terms equal to~$1$;

b) each integer in the interval $[2\ell-S,S]$ is representable as
the sum of a subsequence of~$\alpha$ with length at least~$2\ell-S$.
\end{prop}

\begin{pf}
Part~a) is straightforward. If $\alpha$ contains $x$ terms equal
to~1 then each of the remaining $\ell-x$~terms is at least~2, hence
$S\ge x+2(\ell-x)=2\ell-x$. This implies $x\ge 2\ell-S$. For
part~b), fix $2\ell-S$ ones in~$\alpha$. The remaining
${\ell-(2\ell-S)=S-\ell}$ terms add up to $S-(2\ell-S)=2(S-\ell)$,
so their average is~2. Label these terms $a_1,\dots,a_{S-\ell}$,
assuming that $1\le a_1\le \cdots\le a_{S-\ell}$. Due to this
nondecreasing order, the sequence
$a_1,(a_1+a_2)/2,(a_1+a_2+a_3)/3,\dots$ of arithmetic means is
nondecreasing, hence these means are all at most~2. In other words,
$a_1+\cdots+a_i\le 2i$ for all $i=1,\dots,S-\ell$.

Suppose that $b_1, \dots,b_k$ are \pin/s such that
$b_1+\cdots+b_i\le 2i$ for all $i=1,\dots,k$. Denoting
$S_k=\sum_{i=1}^{k}b_i$, we prove by induction on~$k$ that the
sumset of the sequence $1b_1\dots b_k$ is $\{1,2,\dots,S_k+1\}$.
By {\em sumset} we mean the set of integers representable as a
nonempty subsequence sum. The base $k=1$ is clear. For the
inductive step, let $\Sigma_{k-1}$ and $\Sigma_{k}$ be the sumsets
of $1b_1\dots b_{k-1}$ and $1b_1\dots b_{k-1}b_{k}$, respectively.
Now $\Sigma_{k-1}=\{1,2,\dots,S_{k-1}+1\}$ by the induction
hypothesis, hence
$\Sigma_{k}=\{1,2,\dots,S_{k-1}+1\}\cup\{b_{k},b_{k}+1,\dots,b_{k}+S_{k-1}+1\}$.
Since $b_{k}+S_{k-1}=S_k$, it suffices to check that $b_k\le
S_{k-1}+2$ which is equivalent to $S_k\le 2S_{k-1}+2$. The latter
is true as $2S_{k-1}+2\ge 2(k-1)+2=2k\ge S_k$. The induction is
complete.

Going back to the proof of~b), we infer from the previous
paragraph that the sequence $1a_1\dots a_{S-\ell}$ has sumset
$\{1,2,\dots,2(S-\ell)+1\}$. Take an arbitrary $x\in [2\ell-S,S]$
and set $y=x-(2\ell-S-1)$. Since $1\le y\le 2(S-\ell)+1$, one can
express~$y$ as a nonempty subsequence sum of $1a_1\dots
a_{S-\ell}$. In view of~a), adding $2\ell-S-1$ to both sides of
this representation shows that $x$ equals the sum of at least
$2\ell-S$~terms of the original sequence~$\alpha$. \qed
\end{pf}

\section{The main result}
\label{MainResult}

We are about to characterize all sufficiently long $n$-\zfs/s
in~$\Zn/$. Up to similitude, a sequence of length greater
than~$3n/2-1$ is $n$-\zf/ if and only if it can be divided into two
sequences $\alpha$ and~$\beta$ satisfying $L(\alpha)<n$ and
$L(1-\beta)<n$. (Recall that $L(\omega)$ denotes the sum of the
\lpr/ of the sequence~$\omega$.) There exist sequences of any length
less than $2n-1$ that are ``separable" in the sense just described.
We discuss them before the main theorem in order to indicate that
most of their basic properties do not depend on whether or not the
sequence is ``long."

A couple of technical remarks will be necessary. Let $\alpha$
and~$\beta$ be sequences in~$\Zn/$ such that $L(\alpha)<n$ and
$L(1-\beta)<n$. Because $\v{0}=n$, note that $a\ne 0$ for
$a\in\alpha$ and $b\ne 1$ for $b\in\beta$. We will need the
observations that
\begin{equation}\label{technical}
\v{-b}=n-\v{b}\quad\text{and}\quad \v{1-b}=1+\v{-b}\qquad\text{for
each $b\in\Zn/$, $b\ne 0$.}
\end{equation}
By~(\ref{technical}), for each sequence $\beta$ in~ $\Zn/$ one can
write
\begin{equation}\label{bars}
L(1-\beta)=\sum_{b\in\beta,b\ne 0}\v{-b}+|\beta|.
\end{equation}
In what follows, the empty sequence is assumed to have sum~$0$,
in~$\Z/$ and in~$\Zn/$.

\begin{prop}\label{separable}
Let $n$ and $k$ be integers such that $0<k<n$. Suppose that the
sequences $\alpha$ and~$\beta$ in~$\Zn/$ satisfy the conditions
$|\alpha|+|\beta|=n-1+k$, $L(\alpha)<n$ and $L(1-\beta)<n$. Then:
\begin{enumerate}
\item[a)] The union $\alpha\cup\beta$ is $n$-\zf/.

\item[b)] $k\le |\alpha|<n$, $k\le |\beta|<n$ and $\v{b}-\v{a}\ge k$
for all $a\in\alpha$, $b\in\beta$. In particular $a\ne b$ for all
$a\in\alpha$, $b\in\beta$.

\item[c)] The multiplicities $u$ and $v$ of \kern1pt$1$ and \kern1pt$0$
in~$\alpha\cup\beta$ satisfy
\[
u+v\ge 2k,\quad \max(u,v)\ge k,\quad \min(u,v)\ge 2k-n+1.
\]
The equality $u+v=2k$ is attained if and only if
$\alpha=1^{2p-n+1}2^{n-1-p}$ and $\beta=0^{2q-n+1}(-1)^{n-1-q}$,
for integers $p$ and $q$ such that $(n-1)/2\le p<n$, $(n-1)/2\le
q<n$ and $p+q=n-1+k$. The equality $\max(u,v)=k$  is attained if
and only if $n$ and $k$ have different parity and
$\alpha=1^{k}2^{(n-1-k)/2}$, $\beta=0^{k}(-1)^{(n-1-k)/2}$.

\looseness=-1 \item[d)] For $k\ge (n{-}1)/2$, the highest
multiplicity of a term in~$\alpha\cup\beta$ is $\max(u,v)$.
\end{enumerate}
\end{prop}

\begin{pf}
a) Consider a zero subsequence~$\gamma$ of~$\alpha\cup\beta$. Let
$\gamma$ contain $r$~terms $a_1,\dots ,a_r$ from~$\alpha$,
$s$~nonzero terms $b_1,\dots,b_s$ from~$\beta$, and several zeros,
from~$\beta$ again. Because the sum of $\gamma$ is zero in~$\Zn/$,
the integers $\sum_{i=1}^r\v{a_i}$ and $\sum_{j=1}^s\v{-b_j}$ are
congruent modulo~$n$. Also $0\le \sum_{i=1}^r\v{a_i}\le L(\alpha)<n$
and, by~(\ref{bars}),
\[0\le \sum_{j=1}^s\v{-b_j}\le \sum_{b\in\beta,b\ne
0}\v{-b}=L(1-\beta)-|\beta|<n-|\beta|\le n.
\]
Hence $\sum_{i=1}^r\v{a_i}=\sum_{j=1}^s\v{-b_j}$. Therefore $r\le
\sum_{i=1}^r\v{a_i}=\sum_{j=1}^s\v{-b_j} <n-|\beta|$, implying
$r+|\beta|<n$. Since $|\gamma|\le r+|\beta|$, we infer that
$\alpha\cup\beta$ is $n$-\zf/.

b) The first two inequalities are immediate, because
$|\alpha|+|\beta|=n-1+k$, $|\alpha|\le L(\alpha)<n$ and
$|\beta|=|1-\beta|\le L(1-\beta)<n$. To show that $\v{b}-\v{a}\ge k$
for $a\in\alpha$ and $b\in\beta$, denote
$M=\max_{a\in\alpha}\v{a}+\max_{b\in\beta}\v{1-b}$. Then
\[
2(n-1)\ge L(\alpha)+L(1-\beta)\ge
M+(|\alpha|-1)+(|\beta|-1)=M+n-3+k.
\]
This yields $M\le n+1-k$; thus $\v{a}+\v{1-b}\le n+1-k$ for all
$a\in\alpha$, $b\in\beta$. If $b\ne 0$ then
$\v{1-b}=1+\v{-b}=1+n-\v{b}$ by~(\ref{technical}), so
$\v{a}+\v{1-b}\le n+1-k$ becomes $\v{b}-\v{a}\ge k$. The same
conclusion holds if $b=0$, as then $\v{b}=n$, $\v{1-b}=1$.

c) We have $n-1 \ge L(\alpha)\ge u+2(|\alpha|-u)=2|\alpha|-u$, since
$\v{a}\ge 2$ for $a\ne 1$. Similarly, $\v{1-b}\ge 2$ for $b\ne 0$,
so that $n-1 \ge L(1-\beta)\ge v+2(|\beta|-v)=2|\beta|-v$. Adding up
yields $2(n-1)\ge 2(|\alpha|+|\beta|)-(u+v)=2(n-1+k)-(u+v)$. It
follows that $u+v\ge 2k$, so $\max(u,v)\ge k$. Clearly $\max(u,v)\le
n-1$ by~b), which implies $\min(u,v)\ge 2k-n+1$.

The equality $u+v=2k$ occurs if and only if
$n-1=L(\alpha)=2|\alpha|-u$ and $n-1 =L(1-\beta)=2|\beta|-v$.
These conditions imply $u=2|\alpha|-n+1$, $v=2|\beta|-n+1$; also
$\v{a}=2$ for $a\in\alpha$, $a\ne 1$ and $\v{1-b}=2$ for
$b\in\beta$, $b\ne 0$. In particular $|\alpha|\ge (n-1)/2$,
$|\beta|\ge (n-1)/2$. So setting $p=|\alpha|$, $q=|\beta|$ and
taking~b) into account, we obtain $(n-1)/2\le p<n$, $(n-1)/2\le
q<n$, $p+q=n-1+k$ and $\alpha=1^{2p-n+1}2^{n-1-p}$,
$\beta=0^{2q-n+1}(-1)^{n-1-q}$. The converse is easy to check; we
note only that the last two sequences are well-defined whenever
$(n-1)/2\le p<n$, $(n-1)/2\le q<n$.

If $\max(u,v)=k$ then ${u=v=k}$ in view of $u+v\ge 2k$, so
$u+v=2k$. The conclusions of the last paragraph imply
${\alpha=1^{k}2^{(n-1-k)/2}}$, $\beta=0^{k}(-1)^{(n-1-k)/2}$.
These sequences are well-defined only if $k\not\equiv n\
\text{(mod 2)}$. The converse is clear.

d) We have $u+v\ge 2k$ by~c), so the number of terms different from
$1$ and $0$ in~$\alpha\cup\beta$ is at most $(n-1+k)-2k=n-1-k$. Now
it suffices to observe that $n-1-k\le k$ for $k\ge (n-1)/2$, and
that $\max(u,v)\ge k$ by~c). \qed
\end{pf}

One can see that sequences ``separable" in the sense of
Proposition~\ref{separable} are not just $n$-\zf/ but have an
interesting general structure. This is unexpected at first glance
as $\alpha$ and $\beta$ do not seem to be related in any way.
However, while the properties listed in
Proposition~\ref{separable} are fairly simple to derive, it is
less trivial to establish that each sufficiently long $n$-\zfs/
in~$\Zn/$ is ``separable." The next theorem proves that length
greater than~${3n/2{-}1}$ is enough to guarantee this. Moreover,
shorter $n$-\zfs/s are not necessarily ``separable." These
conclusions form the essential part of the article.

\begin{thm}\label{characterization}
A sequence of length greater than~$3n/2{-}1$ in the cyclic
group~$\Zn/$ does not contain an $n$-term zero subsequence if and
only if it is similar to the union of two sequences $\alpha$ and
$\beta$ in~$\Zn/$ such that
\begin{equation*}
L(\alpha)<n\qquad\text{and}\qquad L(1-\beta)<n.
\end{equation*}
\end{thm}

\begin{pf}
The sufficiency follows from Proposition~\ref{separable}~a). For
the necessity, let $\gamma$ be an $n$-\zfs/ of length greater
than~$3n/2{-}1$ in~$\Zn/$. Translations in~$\Zn/$ do not affect
sums of length~$n$, so one may assume that 0 is a term of~$\gamma$
with maximum multiplicity~$v$. Then Proposition~\ref{gao} shows
that each zero subsequence of~$\gamma$ has length less than~$n$.
In particular $v<n$.

Select a zero subsequence~$\sigma$ of $\gamma$ with nonzero terms
and of maximum length; $\sigma$~may be the empty sequence. This
choice implies that the remaining nonzero terms of~$\gamma$ form a
\zfs/~$\tau$. By the remark above, the lengths of~$\sigma$
and~$\tau$ satisfy $|\sigma|<n-v$ and
$|\tau|>(3n/2{-}1)-(n{-}1)=n/2$. Therefore Theorem~\ref{zero-free}
applies to the \zfs/~$\tau$.

Hence multiplying~$\tau$ by a certain integer coprime to~$n$
yields an equivalent sequence with sum of the \lpr/ less than~$n$.
We multiply by the same integer all remaining terms of~$\gamma$,
which preserves the zero sums of any length. So there is no loss
of generality in assuming that $\gamma=0^v\sigma\tau$, where
$\sigma$ is a zero subsequence of~$\gamma$ with nonzero terms and
of maximum length, and $\tau$ is a \zfs/ of length greater
than~$n/2$ satisfying $L(\tau)<n$.

Let $\sigma=1^wb_1\dots b_q$ where $b_1,\dots,b_q$ are all terms of
$\sigma$ different from~$1$. The following inequality stronger than
$L(\tau)<n$ implies the conclusion directly:
\begin{equation}\label{claim2}
L(\tau)+\sumb<n.
\end{equation}
Indeed, assume~(\ref{claim2}) is true, and let $\alpha=1^w\tau$,
$\beta=0^vb_1\dots b_q$. Then $\gamma=\alpha\cup\beta$; in
addition, $L(\alpha)<n$ and $L(1-\beta)<n$. Firstly, $w\equiv
\sumb\ \text{(mod~$n$)}$, since $\sigma$ has sum zero. Also
$0\le\sumb<n$ by~(\ref{claim2}), and clearly $0\le w<n$. Therefore
$w=\sumb$. So (\ref{claim2}) can be written as $L(\tau)+w<n$,
which is the inequality~$L(\alpha)<n$. Furthermore, we obtain
$|\sigma|=w+q=\sumb+q$, and because $|\sigma|<n-v$, it follows
that $\sumb+q<n-v$. This is the same as $\sum_{b\in\beta,b\ne
0}\v{-b}+|\beta|<n$. By~(\ref{bars}), the latter means that
$L(1-\beta)<n$ .

\looseness=-1 So it suffices to prove~(\ref{claim2}) which is clear
if $\sigma$ is the empty sequence. Hence let $\sigma\ne\emptyset$,
implying $q\ge 1$ ($\sigma$ cannot have only terms equal to~$1$).
For the sake of clarity, denote $|\tau|=\ell>n/2$, $L(\tau)=S<n$ and
$\v{-b_j}=v_j$, $j=1,\dots,q$. Note that $2\ell-S\ge 2\ell-(n-1)\ge
2$ as $\ell>n/2$. Thus Proposition~\ref{general} applies to the
sequence~$\v{\tau}$ of the \lpr/ of~$\tau$. Also, $1\le v_j<n{-}1$
by the choice of $b_1,\dots,b_q$. The proof of~(\ref{claim2}) is
based on the following observation.
\begin{quote}
\em Suppose that $m$~terms $v_{j_1},\dots ,v_{j_m}$ of the
sequence $v_{1}\dots v_{q}$ are such that the integer
$T=n-(v_{j_1}+\cdots+v_{j_m})$ satisfies $1<T\le S$. Then $m\ge
2\ell-S$ if $2\ell-S\le T\le S$ and $m\ge T$ if $1<T<2\ell-S$.
\end{quote}
Indeed, if $T$ represents the congruence class~$t$ modulo~$n$ then
$t=\sum_{i=1}^{m}b_{i_{j}}$. Let $2\ell-S\le T\le S$. By
Proposition~\ref{general}, there is a subsequence~$\omega$
of~$\tau$ with length at least $2\ell-S$ such that
$T=\sum_{c\in\omega}\v{c}$. Hence
$\sum_{c\in\omega}c=t=\sum_{i=1}^{m}b_{i_{j}}$. This implies
$m\ge|\omega|\ge 2\ell-S$ as $m<|\omega|$ would yield a zero
subsequence of $\gamma$ with nonzero terms which is longer than
$\sigma$, obtained upon replacing $b_{j_1}\dots b_{j_m}$
by~$\omega$. Similarly, if $1<T<2\ell-S$ then $T$ can be expressed
as the sum of $T$ terms equal to~$1$ of~$\v{\tau}$. (There are at
least~$2\ell-S$ ones in $\v{\tau}$ by Proposition~\ref{general}.)
Now the same argument as above gives $m\ge T$, by the maximum
choice of~$\sigma$.

It follows from the observation that $n-v_j>S$, $j=1,\dots,q$.
Indeed, if $1<n-v_j\le S$ for some~$j$ then $1\ge 2\ell-S$ or $1\ge
n-v_j$, both of which are not true. Therefore $1\le v_j<n-S$,
$j=1,\dots,q$.

Passing on to the proof of~(\ref{claim2}), suppose that it is
false. Then there are subsequences of $v_1\dots v_q$ whose sum is
at least~$n-S$, for instance $v_1 \dots v_q$ itself. Without loss
of generality, let $v_1\dots v_m$ be such a (nonempty) subsequence
of minimum length~$m$. So $T=n-\sum_{j=1}^{m}v_j\le S$ but
$T+v_j>S$ for all $j=1,\dots,m$. Note that $T>1$ in view of the
previous paragraph, because $v_m<n-S$ yields
$T>S-v_m>S-(n-S)=2S-n\ge 2\ell-n\ge 1$.

Let $2\ell-S\le T\le S$. Then $m\ge 2\ell-S$ by the observation
above. Hence
\[
S+1\le n-\sum_{j=1}^{m-1}v_j \le n-(m-1)\le
n-(2\ell-S-1)=(n-2\ell)+S+1,
\]
implying $n\ge 2\ell$ which is a contradiction.

Next, let $1<T<2\ell-S$. Now the observation gives $m\ge T$.
Recalling that $T+v_j>S$, we have $v_j\ge S+1-T>0$ for
$j=1,\dots,m$, implying
\[
n=T+\sum_{j=1}^mv_j\ge T+m(S+1-T)\ge T+T(S+1-T)=T(S+2-T).
\]
Consider the quadratic function $g(t)=t^2-(S+2)t+n$. We obtained
$g(T)\ge 0$ for some $T\in\{2,\dots,2\ell-S-1\}$. But the maximum
of~$g$ on $\{2,\dots,2\ell-S-1\}$ is $g(2)=n-2S$, and $n-2S\le
n-2\ell<0$. This is a contradiction again; claim~(\ref{claim2})
follows, concluding the main argument. \qed
\end{pf}

Theorem~\ref{characterization} establishes the desired
characterization, in a form hopefully providing general insight
into the structure of $n$-\zfs/s.  On the other hand, the
practically important consequence of the theorem is that each
$n$-\zfs/ of length~${n-1+k}$ in~$\Zn/$, where ${n/2<k<n}$, is
similar to a sequence satisfying the conclusions of
Proposition~\ref{separable}. Both Theorem~\ref{characterization}
and Proposition~\ref{separable} are needed for a really clear
picture of the ``long" $n$-\zfs/s in~$\Zn/$. The next observation
adds one more detail to this picture.

The affine map $x\mapsto 1-x$ interchanges~0 and~1 and transforms
arbitrary sequences $\alpha$ and $\beta$ into $\alpha_1=1-\alpha$
and $\beta_1=1-\beta$, respectively. If the inequalities
$L(\alpha)<n$ and $L(1-\beta)<n$ hold true, they can be written as
$L(1-\alpha_1)<n$ and $L(\beta_1)<n$. So if a sequence $\gamma$ is
similar to $\alpha\cup\beta$, it is also similar to
$\alpha_1\cup\beta_1$, a sequence with all properties from
Proposition~\ref{separable}, in which the multiplicities of~0 and~1
are interchanged. Therefore one can assume additionally that $u\le
v$. For $k\ge (n-1)/2$, Proposition~\ref{separable}~d) then implies
that $0$ is a term of highest multiplicity in $\alpha\cup\beta$.

The conditions $L(\alpha)<n$ and $L(1-\beta)<n$ can be expanded to
obtain an explicit form of the characterization established in
Theorem~\ref{characterization}. Up to certain details, this
explicit description  has the same shape as the one
in~\cite{BialostockiDierkerGrynkiewiczLotspeich} of the $n$-\zfs/s
with length $n-1+k$, for $k$ roughly greater than~$3n/4$. It is
worth noting that the range~$n/2<k<n$ for~$k$ is the natural scope
of such a characterization. There are $n$-\zfs/s of length
$n-1+\lfloor n/2\rfloor$ that do not obey the conclusion of
Theorem~\ref{characterization}.

Here are examples. For an odd $n\ge 9$, consider the sequence
$0^{n-1}2^{(n-5)/2}3^2$. Its length is $n-1+(n-1)/2=n-1+\lfloor
n/2\rfloor$. For an even $n\ge 6$, consider the sequence
$0^{n-1}2^{n/2-1}3$, of length $n-1+n/2=n-1+\lfloor n/2\rfloor$.
Both sequences are $n$-\zf/. Suppose that either of them is
similar to a union $\alpha\cup\beta$ where $\alpha$ and $\beta$
satisfy~$L(\alpha)<n$, $L(1-\beta)<n$. Because $k\ge (n-1)/2$,
$\alpha$ and $\beta$ can be chosen so that $0$ is a term of
highest multiplicity~$v$ in~$\alpha\cup\beta$. Then $v=n-1$, so
$\beta=0^{n-1}$. It follows that $2^{(n-5)/2}3^2$ or $2^{n/2-1}3$
is equivalent to $\alpha$, a sequence satisfying $L(\alpha)<n$.
However, one can check that the latter is not true.

\section{Terms of high multiplicity}
\label{Highmult}

Let $n/2<k<n$, and let $\gamma$ be an $n$-\zfs/ of length $n-1+k$ in
$\Zn/$. It follows from Theorem~\ref{characterization} and
Proposition~\ref{separable} that $\gamma$ contains a term of
multiplicity at least $k$, and two distinct terms of combined
multiplicity at least $2k$. Now we obtain precise forms of these
statements.

Let $\alpha\cup\beta$ be a sequence similar to~$\gamma$ and
satisfying the conclusions of Proposition~\ref{separable}. We may
also assume that $0$ is a term of maximum multiplicity~$v$
in~$\alpha\cup\beta$, as explained in the previous section.

By Proposition~\ref{separable}~c), the equality $v=k$ holds if and
only if $k\not\equiv n\ \text{(mod 2)}$ and $\alpha\cup\beta=0^k1^k2^{(n-1-k)/2}(-1)^{(n-1-k)/2}$. This is the unique sequence
satisfying $v=k$, up to affine maps and rearrangement of terms.

Suppose now that $k\equiv n\ \text{(mod 2)}$. Then $v\ge k+1$ by
Proposition~\ref{separable}~c) again. The equality $v= k+1$ can be
attained, for instance for the  sequence
$0^{k+1}1^{k-1}2^{(n-k)/2}(-1)^{(n-k-2)/2}$ which is well-defined
and $n$-\zf/ by Proposition~\ref{separable}~a) (setting
$\alpha=1^{k-1}2^{(n-k)/2}$, $\beta=0^{k+1}(-1)^{(n-k-2)/2}$).

Thus the following corollary was proved.

\begin{cor}\label{bestmult}
Let $n$ and~$k$ be integers satisfying $n/2<k<n$. Each $n$-\zfs/ of
length $n-1+k$ in~$\Zn/$ contains a term of multiplicity at
least~$k$, if $n$ and $k$ are of different parity, and at
least~$k+1$, if $n$ and $k$ are of the same parity. Both estimates
are best possible.
\end{cor}

The sum of the two highest multiplicities was probably the most
widely explored question concerning $n$-\zfs/s in~$\Zn/$. We are now
in the position to resolve this question completely for each length
$n-1+k$ where $n/2<k<n$. Indeed, the lower bound $u+v\ge 2k$ for
this combined multiplicity follows from above. Now let us take
another look at the examples for the maximum multiplicity of a
single term. In both possible cases, $k\not\equiv n\ \text{(mod 2)}$
and $k\equiv n\ \text{(mod 2)}$, it is easy to see that 0 and 1 are
the two terms with highest combined multiplicity, and the value of
this multiplicity is~$2k$. So we proceed with one more structural
conclusion.

\begin{cor}\label{bestcombinedmult}
Let $n$ and~$k$ be integers satisfying $n/2<k<n$. Each $n$-\zfs/ of
length $n-1+k$ in~$\Zn/$ contains two terms of combined multiplicity
at least~$2k$, and this estimate is best possible.
\end{cor}

Naturally, some well-known results about the structure of the
$n$-\zfs/s with a certain length are now immediate. For example, let
us consider the lengths~$2n{-}2$ (as in~\cite{PetersonYuster},
\cite{BialostockiDierker}) and~$2n{-}3$ (as in~\cite{FloresOrdaz}).
By the discussion above, any $n$-\zfs/ of length~$2n{-}2$ (i.~e.
$n-1+k$ with $k=n{-}1$) is similar to~$0^v1^u$, where $v=n{-}1$ (as
$v\ge k$) and $u=n{-}1$ (as $u+v\ge 2k$). Here we assume $n>2$ to
ensure that $k>n/2$; however the same conclusion holds true for
$n=2$ as well. Similarly, for $n>4$, any $n$-\zfs/ of
length~$2n{-}3$ (i.~e. $n-1+k$ with $k=n{-}2$) is similar
to~$0^v1^u\gamma$, where $v$ is the maximum multiplicity of a term
and $u+v\ge 2k=2n-4$. Since $k=n{-}2\equiv n\ \text{(mod 2)}$,
Corollary~\ref{bestmult} implies $v\ge k{+}1=n{-}1$. Hence $v=n{-}1$
and $u=n{-}2$ or $u=n{-}3$. Now it is easy to infer that each
$n$-\zfs/ of length~$2n{-}3$, $n>4$, is similar to~$0^{n-1}1^{n-2}$
or~$0^{n-1}1^{n-3}2$. For $n=3,4$, this conclusion can be checked
directly. For~$n=2$, the only $n$-\zfs/ of length $2n-3=1$ is
similar to the one-term sequence~$0$.

\section{The $g(n,k)$ function}
\label{gnk}

For \pin/s $n$ and $k$, $k\le n$, let $g(n,k)\ge k$ be the least
integer such that each sequence in~$\Zn/$ with at least $k$ distinct
terms and length $g(n,k)$ contains an $n$-term zero sum. The
function $g(n,k)$ was introduced by Bialostocki and Lotspeich
in~\cite{BialostockiLotspeich}. The structural results about
$n$-\zfs/s of lengths~$2n{-}2$ and~$2n{-}3$ (such as mentioned after
Corollary~\ref{bestcombinedmult}) imply $g(n,2)=2n{-}1$ for $n\ge 2$
and $g(n,3)=2n{-}2$ for~$n\ge 4$. It is easy to see that $g(3,3)=3$.

The values of $g(n,k)$ for $4\le k\le \sqrt{n+4}+1$ were found
in~\cite{BialostockiDierkerGrynkiewiczLotspeich}: If $k\ge 4$ is
even and $n\ge k^2-2k-4$, or if $k\ge 5$ is odd and $n\ge
k^2-2k-3$, then
\begin{equation*}
g(n,k)=2n-1-\intpart.
\end{equation*}
For the lower bound, the following examples were used. In the case
of an even~$k\ge 2$, consider the sequence
\[
-\left(
\frac{k-2}{2}\right)\dots(-1)(0)^{n-(k^2+2k)/8}(1)^{n-(k^2+2k)/8}(2)\dots\left(\frac
k2\right);
\]
if $k\ge 3$ is odd, consider the sequence
\[
-\left(
\frac{k-3}{2}\right)\dots(-1)(0)^{n-(k^2-1)/8}(1)^{n-(k^2+4k+3)/8}(2)\dots\left(\frac
{k+1}{2}\right).
\]
These examples are valid under the weaker restrictions $n\ge
(k^2+2k)/8+1$ when $k$ is even, and $n\ge (k^2+4k+3)/8+1$ when $k$
is odd. The multiplicities of~0 and~1 in both sequences are
positive integers. By Proposition~\ref{separable}~a), both
sequences are $n$-\zf/, and each one contains $k$ distinct terms.

Here we prove that the function $g(n,k)$ obeys the same formula as
above under the weaker constraints $4\le k\le\sqrt{2n-1}+1$. In this
range the examples above still provide the lower bound $g(n,k)\ge
2n-1-\lfloor((k{-}1)/2)^2\rfloor$.

\begin{thm}
\label{brakemeier} Let $n\ge k$ be integers such that ${4\le
k\le\sqrt{2n-1}+1}$. Then
\begin{equation*}
g(n,k)=2n-1-\intpart.
\end{equation*}
\end{thm}

\begin{pf}
As already mentioned, we need to prove only the upper bound. The
condition $k\le\sqrt{2n-1}+1$ is equivalent to
$n-\lfloor((k{-}1)/2)^2\rfloor>n/2$. Also $k\ge 4$, so the integer
$\ell=n-\lfloor((k{-}1)/2)^2\rfloor$ satisfies $n/2<\ell<n$.
Consider any $n$-\zfs/ $\gamma$ of length $n-1+\ell$ in~$\Zn/$. It
suffices to prove that the number of distinct terms in~$\gamma$ is
less than $k$; then the definition of $g(n,k)$ implies $g(n,k)\le
n-1+\ell=2n-1-\lfloor((k{-}1)/2)^2\rfloor$.

Let $\alpha\cup\beta$ be a sequence similar to $\gamma$ , where
$\alpha$ and~$\beta$ satisfy the conditions in
Proposition~\ref{separable}, with $k$ replaced by~$\ell$. Let there
be $x$ distinct terms in~$\alpha$ and $y$ distinct terms in~$\beta$.
Then Proposition~\ref{separable}~b) shows that the number of
distinct terms in $\gamma$ is $z=x+y$. The sum $L(\alpha)$ does not
increase upon replacing $x$ distinct summands in it by the least
possible values $1,2,\dots,x$, and all remaining summands by~$1$.
Therefore
\[
1+2+\cdots+x+(|\alpha|-x)\le L(\alpha)\le n-1,
\]
which gives $(x^2-x)/2\le n-1-|\alpha|$. Likewise, noticing that
there are $y$ distinct terms in $1-\beta$, we obtain $(y^2-y)/2\le
n-1-|\beta|$. Hence
\[
\frac12(x^2-x)+\frac12(y^2-y)\le
2(n-1)-\left(|\alpha|+|\beta|\right).
\]
Because $|\alpha|+|\beta|=n-1+\ell$, the right-hand side expression
is equal to $n-1-\ell=\lfloor((k{-}1)/2)^2\rfloor-1$. On the other
hand,
\[
\frac12(x^2-x)+\frac12(y^2-y)\ge
\left(\frac{x+y}{2}\right)^2-\frac{x+y}{2}=\left(\frac{z}{2}\right)^2-\frac{z}{2}
=\left(\frac{z-1}{2}\right)^2-\frac 14.
\]
Thus $((z{-}1)/2)^2-1/4\le \lfloor((k{-}1)/2)^2\rfloor -1$ which
implies the desired $z<k$ and completes the proof.\qed
\end{pf}

\end{document}